\documentclass[11pt]{article}

\usepackage{amssymb,amscd,array}

\let\ssection=\section
\renewcommand{\section}{\setcounter{equation}{0}\ssection}

\newcommand{\bbR}{\mathbb{R}}
\newcommand{\bbRP}{\mathbb{RP}}

\newcommand{\Diff}{\mathrm{Diff}}

\newcommand{\cF}{{\mathcal{F}}}
\newcommand{\cD}{{\mathcal{D}}}
\newcommand{\cT}{{\mathcal{T}}}

\newcommand{\Hom}{\mathrm{Hom}}

\newcommand{\PSL}{\mathrm{PSL}}
\newcommand{\SL}{\mathrm{SL}}
\newcommand{\Sl}{\mathrm{sl}}

\newcommand{\Vect}{\mathrm{Vect}}

\newcommand{\cqfd}{\hspace*{\fill}\rule{3mm}{3mm}}

\begin{document}



\def\d{\delta}
\def\g{\gamma}
\def\om{\omega}
\def\r{\rho}
\def\a{\alpha}
\def\b{\beta}
\def\s{\sigma}
\def\vfi{\varphi}
\def\l{\lambda}
\def\m{\mu}
\def\implies{\Rightarrow}

\oddsidemargin .1truein
\newtheorem{thm}{Theorem}[section]
\newtheorem{lem}[thm]{Lemma}
\newtheorem{cor}[thm]{Corollary}
\newtheorem{pro}[thm]{Proposition}
\newtheorem{ex}[thm]{Example}
\newtheorem{rmk}[thm]{Remark}
\newtheorem{defi}[thm]{Definition}


\title{Schwarzian derivative related to modules
of differential operators on a locally projective manifold}

\author{
S. Bouarroudj
\thanks{mailto:sofbou@cpt.univ-mrs.fr}
\and
V.Yu. Ovsienko\thanks{mailto:ovsienko@cpt.univ-mrs.fr
\hfill\break
Centre de Physique Th\'eorique,
CPT-CNRS, Luminy Case 907,
F--13288 Marseille, Cedex 9, FRANCE.
}
}
\date{}
\maketitle 

\begin{abstract}
We introduce a 1-cocycle on the group of diffeomorphisms $\Diff(M)$
of a smooth manifold $M$ endowed with a projective connection.
This cocycle represents a nontrivial cohomology class of $\Diff(M)$
related to the $\Diff(M)$-modules of second order linear differential 
operators on $M$. In the one-dimensional case, this cocycle coincides 
with the Schwarzian derivative, while, in the multi-dimensional case, it
represents its natural and new generalization.
This work is a continuation of \cite{bo} where the same problems
have been treated in one-dimensional case.
\end{abstract}

\section{Introduction}

{\bf 1.1 The classical Schwarzian derivative.} 
Consider the group $\Diff(S^1)$ of diffeomorphisms of the circle preserving 
its orientation. Identifying $S^1$ with $\bbRP^1$,
fix an affine parameter $x$ on $S^1$ such that the natural $\PSL(2,\bbR)$-action is
given by the linear-fractional transformations:
\begin{equation}
x\to\frac{ax+b}{cx+d}, 
\quad\hbox{where}\quad
\left(
\begin{array}{cc}
a&b\\
c&d
\end{array}
\right)
\in\SL(2,\bbR).
\label{Moebius}
\end{equation}
The classical Schwarzian derivative is then given by:
\begin{equation}
S(f)=
\left(\frac{f'''(x)}{f'(x)}-
\frac{3}{2}\left(\frac{f''(x)}{f'(x)}\right)^2\right)(dx)^2,
\label{Schwarzian}
\end{equation}
where $f\in \Diff(S^1)$.


\medskip
\noindent
{\bf 1.2 The Schwarzian derivative as a 1-cocycle.}
It is well known that the Schwarzian derivative can be intrinsically
defined as the \textit{unique 1-cocycle} on $\Diff(S^1)$ with values in the space of 
quadratic differentials on $S^1$, \textit{equivariant with respect to the 
M\"obius group} $\PSL(2,\bbR)\subset\Diff(S^1)$, cf.~\cite{bot,kir}.
That means, the map (\ref{Schwarzian}) satisfies the following two
conditions:
\begin{equation}
S(f\circ g)=
g^*S(f)+S(g),
\label{Cocycle}
\end{equation}
where $f^*$ is the natural $\Diff(S^1)$-action on the space of quadratic differentials
and
\begin{equation}
S(f)=
S(g),\qquad g(x)=\frac{af(x)+b}{cf(x)+d}.
\label{Invariant}
\end{equation}
Moreover, the Schwarzian derivative is characterized by 
(\ref{Cocycle}) and (\ref{Invariant}).


\goodbreak
\noindent
{\bf 1.3 Relation to the module of second order differential operators.}
The Schwarzian derivative appeared in the classical literature in closed
relation with differential operators. More precisely, consider the space of
Sturm-Liouville operators: $A_u=-2\left(\frac{d}{dx}\right)^2+u(x)$, where
$u(x)\in C^{\infty}(S^1)$, the action of $\Diff(S^1)$ on this space
is given by $f(A_u)=A_v$ with
\begin{equation}
v=
u\circ f^{-1}\cdot({f^{-1}}^{\prime}){}^2+
\frac{f'''(x)}{f'(x)}-
\frac{3}{2}\left(\frac{f''(x)}{f'(x)}\right)^2
\label{fu}
\end{equation}
(see e.g. \cite{wil,ca}). 

It, therefore, seems to be clear that the natural
approach to understanding of multi-dimensional analogues of the Schwarzian
derivative should be based on the relation with modules of differential
operators.


\medskip
\noindent
{\bf 1.4 The contents of this paper.}
In this paper we introduce a multi-dimensional analogue of the Schwarzian
derivative related to the $\Diff(M)$-modules of differential operators on
$M$.

Following \cite{do} and \cite{lo}, the module of differential operators 
$\cD_{\l,\m}$ will be viewed as a \textit{deformation} of the module
of symmetric contravariant tensor fields on $M$. This approach leads to
$\Diff(M)$-cohomology first evoked in \cite{do}.
The corresponding cohomology of the Lie algebra of vector fields $\Vect(M)$ 
has been calculated in \cite{lo} for a manifold $M$ endowed with a
flat projective structure. We use these results to determine 
the projectively equivariant cohomology of $\Diff(M)$ arising in this context.

Note that multi-dimensional analogous of the Schwarzian
derivative is a subject
already considered in the literature. We will refer
\cite{al,kh,mm,os,ovs1,tab,ret}
for various versions of multi-dimensional
Schwarzians in projective, conformal,
symplectic and non-commutative geometry.

\section{Projective connections}

Let $M$ be a smooth (or complex) manifold of dimension $n$. 
There exists a notion of projective connection on $M$,
due to E. Cartan.
Let us recall here the simplest (and naive) way to define
a projective connection as an equivalence class of standard
(affine) connections.

\subsection{Symbols of projective connections}

{\bf Definition.}
A \textit{projective connection} on $M$ is the class of affine connections
corresponding to the same expressions
\begin{equation}
\Pi_{ij}^k=
\Gamma_{ij}^k-\frac{1}{n+1}\left(\delta^k_i\Gamma_{jl}^l+
\delta^k_j\Gamma_{il}^l\right),
\label{Crys}
\end{equation} 
where $\Gamma_{ij}^k$ are the Christoffel symbols and we have assumed a summation
over repeated indices.

\medskip

The symbols (\ref{Crys}) 
naturally appear if one consideres projective connections
as a particular case of so-called Cartan normal connection, 
see \cite{kn}. 

\goodbreak
\noindent
{\bf Remarks.}

(a) The definition is correct (i.e. does not depend on the choice of 
local coordinates on $M$). 

(b) The formula (\ref{Crys}) defines a natural projection to
the space of trace-less (2,1)-tensors, one has: $\Pi_{ik}^k=0$.

\subsection{Flat projective connections and projective structures}

A manifold $M$ is said to be locally projective (or endowed with
a \textit{flat projective structure}) if there exists an atlas on $M$ with
linear-fractional coordinate changes~:
\begin{equation}
x^i=
\frac{a^i_jx^j+b^i}{c_jx^j+d}.
\label{LinFrac}
\end{equation} 

A projective connection on $M$ is called \textit{flat} if
in a neighborhood of each point, there exists a local 
coordinate system $(x^1,\ldots,x^n)$ such that the symbols
$\Pi_{ij}^k$ are identically zero (see \cite{kn} for a geometric
definition).
Every flat projective connection defines a projective structure
on $M$.

\subsection{A projectively invariant 1-cocycle on $\Diff(M)$}
\label{FirstCocycle}

A common way of producing nontrivial cocycles on $\Diff(M)$ using
affine connections on $M$ is as follows. The map:
$(f^*\Gamma)_{ij}^k-\Gamma_{ij}^k$ is a 1-cocycle on $\Diff(M)$ with values 
in the space of symmetric $(2,1)$-tensor fields.
It is, therefore, clear that a projective connection on $M$ leads to
the following 1-cocycle on $\Diff(M)$:
\begin{equation}
\ell(f)=
\left(
(f^*\Pi)_{ij}^k-\Pi_{ij}^k
\right)
dx^i\otimes dx^j\otimes\frac{\partial}{\partial x^k}
\label{ell}
\end{equation}
vanishing on (locally) projective diffeomorphisms.

\medskip
\noindent
{\bf Remarks.}

(a)
The expression (\ref{ell}) is well defined (does not depend on the
choice of local coordinates). This follows from a well-known fact, 
that the difference of two (projective) connections defines a 
$(2,1)$-tensor field.

(b)
Already the formula (\ref{ell}) implies that the map $f\mapsto\ell(f)$ is,
indeed, a 1-cocycle, that is, it satisfies the relation
$\ell(f\circ g)=g^*\ell(f)+\ell(g)$. 

(c)
It is clear that the cocycle $\ell$ is nontrivial (cf. \cite{lo}),
otherwise it
would depend only on the first jet of the diffeomorphism $f$. 
Note that the formula (\ref{ell})
looks as a coboundary, however, the
symbols $\Pi_{ij}^k$ do not transform as components of
a $(2,1)$-tensor field (but as symbols of a projective connection).

\medskip
\noindent
{\bf Example.} In the case of a smooth manifold endowed with a
\textit{flat} projective connection, (with
symbols (\ref{Crys}) identically zero) or, equivalently, with
a projective structure, the cocycle (\ref{ell}) 
obviously takes the form:
\begin{equation}
\ell(f,x)=
\left(
\frac{\partial^2f^l}{\partial x^i\partial x^j}
\frac{\partial x^k}{\partial f^l}-
\frac{1}{n+1}\left(\d^k_j\frac{\partial \log J_f}{\partial x^i}
+
\d^k_i\frac{\partial \log J_f}{\partial x^j}
\right)
\right)
dx^i\otimes dx^j\otimes\frac{\partial}{\partial x^k}
\label{barS}
\end{equation}
where $f(x^1,\ldots,x^n)=\left(f^1(x),\ldots,f^n(x)\right)$ and
$J_f=\det\left(\frac{\partial f^i}{\partial x^j}\right)$
is the Jacobian.
This expression is globally defined and
vanishes if $f$ is given (in the local coordinates of the
projective structure) as a linear-fractional transformation
(\ref{LinFrac}).

\medskip
The cocycle (\ref{ell},\ref{barS}) was introduced in 
\cite{tab,mm} as 
a multi-dimensional projective analogue of the Schwarzian derivative.
However, in contradistinction with the Schwarzian derivative 
(\ref{Schwarzian}), 
this map (\ref{barS}) depends on the second-order
jets of diffeomorphisms. 
Moreover, in the one-dimensional case
$(n=1)$, the expression (\ref{ell},\ref{barS}) is identically zero.

\section{Introducing the Schwarzian derivative}

Assume that $\dim M\geq2$.
Let ${\cal S}^k(M)$ (or ${\cal S}^k$ for short)
be the space of $k$-th order symmetric
contravariant  tensor fields on $M$.

\subsection{Operator symbols of a projective connection}

For an arbitrary system of local coordinates fix the following
linear differential operator
$T:{\cal S}^2\to C^{\infty}(M)$ given for every $a\in{\cal S}^2$
by $T(a)=T_{ij}(a^{ij})$ with
\begin{equation}
T_{ij}=
\Pi^k_{ij}\frac{\partial}{\partial x^k}-
\frac{2}{n-1}\,
\left(\frac{\partial\Pi^k_{ij}}{\partial x^k}-
\frac{n+1}{2}\,\Pi^k_{il}\Pi^l_{kj}\right),
\label{MultiSchwar}
\end{equation}
where $\Pi^k_{ij}$ are the symbols of a projective connection
(\ref{Crys}) on $M$.

It is clear that the differential operator (\ref{MultiSchwar}) is
not intrinsically defined, indeed, already its principal symbol,
$\Pi_{ij}^k$, is not a tensor field. In the same spirit that the
difference of two projective connections 
$\widetilde\Pi_{ij}^k-\Pi_{ij}^k$ is a well-defined tensor field,
we have the following

\begin{thm}
Given arbitrary projective connections $\widetilde\Pi_{ij}^k$
and $\Pi_{ij}^k$, the difference
\begin{equation}
{\cal T}=
\widetilde T-T
\label{TheSchwar}
\end{equation}
is a linear differential operator from ${\cal S}^2$ to
$C^{\infty}(M)$ well defined
(globally) on $M$ (i.e., it does not depend on the choice of
local coordinates).
\label{definit}
\end{thm}

\noindent
\textit{Proof.}
To prove that the expression (\ref{TheSchwar}) is, indeed a
well-defined differential operator from 
${\cal S}^2$ into $C^{\infty}(M)$, we need an explicit
formula of coordinate transformation for such kind of
operators.

\begin{lem}
The coefficients of a first-order linear differential
operator $A:{\cal S}^2\to C^{\infty}(M)$
$
A(a)=
\left(t_{ij}^k\partial_k+
u_{ij}\right)a^{ij}
$
transform under coordinate changes as follows:
\begin{eqnarray}
t_{ij}^k(y)&=&
t_{ab}^c(x)\frac{\partial x^a}{\partial y^i}
\frac{\partial x^b}{\partial y^j}
\frac{\partial y^k}{\partial x^c}
\label{transform1}\\
u_{ij}(y)&=&
u_{ab}(x)\frac{\partial x^a}{\partial y^i}
\frac{\partial x^b}{\partial y^j}-
2t_{ab}^c(x)\,
\frac{\partial^2 y^k}{\partial x^c\partial x^l}\,
\frac{\partial x^a}{\partial y^k}
\frac{\partial x^b}{\partial y^{(i}}
\frac{\partial x^l}{\partial y^{j)}}
\label{transform}
\end{eqnarray}
where round brackets mean symmetrization.
\label{lemma}
\end{lem}
\noindent
\textit{Proof of the lemma:}
straightforward.
\cqfd

\goodbreak

Consider the following expression~:
$$
{\cal T}(\alpha,\beta)_{ij}=
\left(\widetilde\Pi_{ij}^k-\Pi_{ij}^k\right)\partial_k+
\alpha\partial_k\left(\widetilde\Pi_{ij}^k-\Pi_{ij}^k\right)+
\beta
\left(
\widetilde\Pi_{li}^k\widetilde\Pi_{jk}^l-
\Pi_{li}^k\Pi_{jk}^l
\right)
$$
From the definition 
(\ref{MultiSchwar},\ref{TheSchwar})
for
\begin{equation}
\alpha=-\frac{2}{n-1},\qquad
\beta=\frac{n+1}{n-1},
\label{solution}
\end{equation}
one gets ${\cal T}(\alpha,\beta)={\cal T}$.

Now, it follows immediately
from the fact that $\widetilde\Pi_{ij}^k-\Pi_{ij}^k$ 
is a well-defined $(2,1)$-tensor
field on $M$,
that the condition (\ref{transform1}) for the principal
symbol of ${\cal T}(\alpha,\beta)$
is satisfied.

The transformation law for the symbols of a
projective connection reads:
$$
\Pi_{ij}^k(y)=
\Pi_{ab}^c(x)\,
\frac{\partial x^a}{\partial y^i}
\frac{\partial x^b}{\partial y^j}
\frac{\partial y^k}{\partial x^c}+
\ell(y,x),
$$
where $\ell(y,x)$ is given by (\ref{barS}).
Let $u(\alpha,\beta)_{ij}$
be the zero-order term in ${\cal T}(\alpha,\beta)_{ij}$,
one readily gets:
\begin{eqnarray*}
u(\alpha,\beta)(y)_{ij}&=&
u(\alpha,\beta)(x)_{ab}
\frac{\partial x^a}{\partial y^i}
\frac{\partial x^b}{\partial y^j}\\
&&-
2(\alpha+\beta)\,
\left(\widetilde{\Pi}_{ab}^c(x)-\Pi_{ab}^c(x)\right)\,
\frac{\partial^2 y^k}{\partial x^c\partial x^l}\,
\frac{\partial x^a}{\partial y^k}
\frac{\partial x^b}{\partial y^{(i}}
\frac{\partial x^l}{\partial y^{j)}}\\
&&+
(\alpha+\frac{2\beta}{n+1})
\left(\widetilde{\Pi}_{ab}^c(x)-\Pi_{ab}^c(x)\right)\,
\frac{\partial \log J_y}{\partial x^c}\,
\frac{\partial x^a}{\partial y^i}
\frac{\partial x^b}{\partial y^j}\,.
\end{eqnarray*}
The transformation law (\ref{transform})
for $u(\alpha,\beta)_{ij}$
is satisfied if and only if
$\alpha$ and $\beta$ are given by (\ref{solution}).
Theorem \ref{definit} is proven.
\cqfd

\medskip

We call $T_{ij}$ given by (\ref{MultiSchwar}) the
\textit{operator symbols} of a projective connection.
This notion is the main tool of this paper.

\medskip
\noindent
\textbf{Remark}.
The scalar term of (\ref{MultiSchwar}) looks similar to the
symbols 
$
\Pi_{ij}=-\partial\Pi^k_{ij}/\partial x^k+
\Pi^k_{il}\Pi^l_{kj},
$
which together with $\Pi^k_{ij}$ characterise
the normal Cartan projective connection (see~\cite{kn}).
We will show that the operator symbols $T_{ij}$,
and not the symbols of the normal projective connection,
lead to a natural notion of multi-dimensional Scharzian derivative.
However, the geometric meaning of (\ref{MultiSchwar}) is still
mysterious for us.

\subsection{The main definition}

Consider a manifold $M$ endowed with a projective connection.
The expression
\begin{equation}
S(f)=
f^*(T)-T,
\label{TheSchwarz}
\end{equation}
where $T$ is the (locally defined) operator (\ref{MultiSchwar}),
is a linear differential operator well defined
(globally) on $M$.

\goodbreak

\begin{pro}
the map
$f\mapsto S(f)$ is a nontrivial 1-cocycle on $\Diff(M)$ with values 
in $\Hom({\cal S}^2,C^{\infty}(M))$.
\label{definition}
\end{pro}

\noindent
\textit{Proof}.
The cocycle property for $S(f)$
follows directly from the definition (\ref{TheSchwarz}).
This cocycle is not a coboundary.
Indeed, every coboundary $\mathrm{d}B$
on $\Diff(M)$ with values in the space 
$\Hom({\cal S}^2,C^{\infty}(M))$
is of the form $B(f)(a)=f^*(B)-B$,
where $B\in\Hom({\cal S}^2,C^{\infty}(M))$.
Since $S(f)$ is a first-order differential operator,
the coboundary condition $S=\mathrm{d}B$ would imply
that $B$ is also a first-order differential operator
and so, $\mathrm{d}B$ depends at most on the second jet of $f$.
But, $S(f)$ depends on the third jet of $f$.
This contradiction proves that the cocycle
(\ref{TheSchwarz}) is nontrivial.
\cqfd

\medskip

The cocycle (\ref{TheSchwarz}) will be called the
\textit{projectively equivariant Schwarzian derivative}.
It is clear that the kernel of $S$ is precisely the
subgroup of $\Diff(M)$ preserving the projective connection.

\medskip
\noindent
{\bf Example.}
In the projectively flat case, $\Pi_{ij}^k\equiv0$,
the cocycle (\ref{TheSchwarz})
takes the form:
\begin{equation}
S(f)_{ij}=
\ell(f)^k_{ij}\frac{\partial}{\partial x^k}-
\frac{2}{n-1}\,\frac{\partial}{\partial x^k}\left(\ell(f)^k_{ij}\right)+
\frac{n+1}{n-1}\,\ell(f)^k_{im}\ell(f)^m_{kj},
\label{MultiSchwar2}
\end{equation}
where $\ell(f)_{ij}^k$ are the components of the cocycle
(\ref{ell}) with
values in symmetric (2,1)-tensor fields.
The cocycle (\ref{MultiSchwar2}) vanishes if and only if $f$ is a
linear-fractional transformation.

It is easy to compute this expression in local coordinates:~:
\begin{equation}
S(f)_{ij}=
\ell(f)^k_{ij}\frac{\partial}{\partial x^k}+
\frac{\partial^3f^k}{\partial x^i\partial x^j\partial x^l}\,
\frac{\partial x^l}{\partial f^k}
-
\frac{n+3}{n+1}\,
\frac{\partial^2 J_f}{\partial x^i\partial x^j}\,J_f^{-1}
+
\frac{n+2}{n+1}\,
\frac{\partial J_f}{\partial x^i}
\frac{\partial J_f}{\partial x^j}\,J_f^{-2}\,.
\label{coord}
\end{equation}
To obtain this formula from (\ref{MultiSchwar2}),
one uses the relation:
$$
\frac{\partial^3f^k}{\partial x^i\partial x^j\partial x^l}\,
\frac{\partial x^l}{\partial f^k}-
\frac{\partial^2f^k}{\partial x^i\partial x^m}\,
\frac{\partial^2f^l}{\partial x^j\partial x^s}\,
\frac{\partial x^m}{\partial f^l}
\frac{\partial x^s}{\partial f^k}=
\frac{\partial^2 J_f}{\partial x^i\partial x^j}\,J_f^{-1}+
\frac{\partial J_f}{\partial x^i}
\frac{\partial J_f}{\partial x^j}\,J_f^{-2}\,.
$$
We observe that, in the one-dimensional case $(n=1)$,
the expression (\ref{coord}) is precisely
$-S(f)$, where $S$ is the classical Schwarzian derivative.
(Recall that in this case $\ell(f)\equiv0$.)

\medskip
\noindent
{\bf Remarks.}

(a)
The infinitesimal analogue of the cocycle (\ref{MultiSchwar2})
has been introduced in \cite{lo}.

(b)
We will show in Section \ref{SecComput}, that the analogue of 
the operator (\ref{TheSchwarz}) in
the one-dimensional case, is, in fact, the operator of multiplication by the
Schwarzian derivative.

\subsection{A remark on the projectively equivariant
cohomology}

Consider the standard $\Sl(n+1,\bbR)$-action on $\bbR^n$
(by infinitesimal projective transformations).
The first group of differential cohomology of $\Vect(\bbR^n)$,
vanishing on the subalgebra $\Sl(n+1,\bbR)$,
with coefficients in the space 
$\cD({\cal S}^k,{\cal S}^\ell)$ of
linear differential operators from ${\cal S}^k$
to ${\cal S}^\ell$, was calculated in \cite{lo}.
For $n\geq2$ the result is as follows:
$$
H^1(\Vect(\bbR^n),\Sl(n+1,\bbR);
\cD({\cal S}^k,{\cal S}^\ell))=
\left\{
\begin{array}{rl}
\bbR,& k-\ell=2,\\
\bbR,& k-\ell=1,\ell\not=0,\\
0,&\hbox{otherwise}
\end{array}
\right.
$$
The cocycle (\ref{MultiSchwar2}) is, in fact, corresponds to the
nontrivial cohomology class in the case $k=2,\ell=0$ integrated
to the group $\Diff(\bbR^n)$, while the
nontrivial cohomology class in the case $k-\ell=1$
is given by the operator of contraction with the tesor field
(\ref{barS}). 

For any locally projective manifold $M$ it follows that
the cocycle (\ref{TheSchwarz}) generates the unique nontrivial
class of the cohomology of $\Diff(M)$ with coefficients in
$\cD({\cal S}^2,C^\infty(M))$,
vanishing on the (pseudo)group of 
(locally defined) projective transformations.
The same fact is true for the cocycle (\ref{ell}).

\section{Relation to the modules of differential operators}

Consider, for simplicity, a smooth oriented manifold $M$.
Denote $\cD(M)$ the space of scalar linear differential operators
$A:C^{\infty}(M)\to C^{\infty}(M)$. There exists a two-parameter 
family of $\Diff(M)$-module structures on $\cD(M)$. To define it,
one identifies the arguments of differential operators with tensor 
densities on $M$ of degree $\l$ and their values with tensor densities on $M$ 
of degree $\m$.

\subsection{Differential operators acting on tensor densities}

Consider the the space $\cF_\l$ of 
\textit{tensor densities} on $M$, that mean, of sections of the line bundle
$({\Lambda^n}T^*M)^\l$. It is clear that $\cF_\l$ is naturally a 
$\Diff(M)$-module.

Since $M$ is oriented, $\cF_\l$ can be identified with
$C^{\infty}(M)$ as a vector space. The $\Diff(M)$-module structures are,
however, different.

\medskip
\noindent
\textbf{Definition.}
We consider the differential operators acting on tensor
densities, namely, 
\begin{equation}
A:\cF_\l\to\cF_\m.
\label{Conv}
\end{equation}
The $\Diff(M)$-action on $\cD(M)$, depending on
two parameters $\l$ and $\m$,
is defined by the usual formula:
\begin{equation}
f_{\l,\m}(A)={f^*}^{-1}\circ A\circ f^*,
\label{Oaction}
\end{equation}
where $f^*$ is the natural $\Diff(M)$-action on $\cF_\l$.

\medskip
\noindent
\textbf{Notation.}
The $\Diff(M)$-module of
differential operators on $M$ with the action (\ref{Oaction}) is denoted
$\cD_{\l,\m}$. For every $k$, the space of differential operators of order
$\leq k$ is a $\Diff(M)$-submodule of $\cD_{\l,\m}$,
denoted $\cD_{\l,\m}^k$.

In this paper we will only deal with the special case $\l=\m$ and use the
notation $\cD_\l$ for $\cD_{\l,\l}$ and $f_\l$ for $f_{\l,\l}$.

\medskip

The modules $\cD_{\l,\m}$ have already been considered in classical
works (see \cite{wil}) and systematically studied in a series of recent papers
(see \cite{do,lmt,lo,bo,do1} and references therein). 

\subsection{Projectively equivariant symbol map}

From now on, we suppose that the manifold $M$ is endowed with
a projective structure.
It was shown in \cite{lo} that there exists a 
(unique up to normalization)
\textit{projectively equivariant symbol map}, that is, a linear 
bijection $\s_\l$
identifying the space $\cD(M)$ with the space of 
symmetric contravariant tensor fields on $M$.

Let us give here the explicit formula
of $\s_\l$ in the case of second order differential operators.
In coordinates of the projective structure, $\s_\l$
associates to a differential operator
\begin{equation}
A=
a_2^{ij}\frac{\partial}{\partial x^i}\frac{\partial}{\partial x^j}+
a_1^i\frac{\partial}{\partial x^i}+
a_0,
\label{Oper}
\end{equation}
where $a_\ell^{{i_1}\ldots{i_\ell}}\in{}C^{\infty}(M)$
with $\ell=0,1,2$,
the tensor field:
\begin{equation}
\sigma_\l(A)=
\bar a_2^{ij}\partial_i\otimes\partial_j
+
\bar a_1^i\partial_i
+\bar a_0,
\label{Tensor}
\end{equation}
and is given by
\begin{equation}
\matrix{
\bar a_2^{ij}
=&
a_2 ^{ij}\hfill \cr\noalign{\smallskip}
\bar a_1^i
=& \displaystyle
a_1^i- 2\frac{(n+1)\lambda+1}{n+3}\,
\frac{\partial a_2^{ij}}{\partial x^j}\hfill
\cr\noalign{\smallskip}
\bar a_0
=& \displaystyle
a_0 - \lambda\,\frac{\partial a_1^i}{\partial x^i} +
\lambda\frac{(n+1)\lambda+1}{n+2}\,
\frac{\partial^2 a_2^{ij}}{\partial x^i\partial x^j} \hfill \cr
}
\label{Symbol}
\end{equation}

The main property of the symbol map $\s_\l$ is that it
commutes with (locally defined) $\SL(n+1,\bbR)$-action.
In other words, the formula (\ref{Symbol}) does not change
under linear-fractional coordinate changes (\ref{LinFrac}).

\subsection{$\Diff(M)$-module of second order 
differential operators}\label{SecComput}

In this section we will compute the $\Diff(M)$-action
$f_\l$ given by (\ref{Oaction}) with $\l=\m$
on the space $\cD^2_\l$ 
(of second order differential operators (\ref{Oper})
acting on $\l$-densities).

Let us give here the explicit formula
of $\Diff(M)$-action in terms of the
projectively invariant symbol $\s^\l$.
Namely, we are looking for the operator 
$\bar f_\l=\s_\l\circ f_\l\circ(\s_\l)^{-1}$
(such that the diagram below is commutative):
\begin{equation}
\begin{CD}
\cD_\l^2 
@> f_\l >>
\cD_\l^2 
\strut\\ 
@V{\s_\l}VV
@VV{\s_\l}V \strut\\
{\cal S}^2\oplus{\cal S}^1\oplus{\cal S}^0 @> \bar f_\l >> 
{\cal S}^2\oplus{\cal S}^1\oplus{\cal S}^0  \strut
\end{CD}
\label{TheDiagram}
\end{equation}
where ${\cal S}^2\oplus{\cal S}^1\oplus{\cal S}^0$
is the space of second order
contravariant tensor fields (\ref{Tensor}) on $M$.

The following statement, whose proof is straightforward,
shows how the cocycles (\ref{ell}) and (\ref{TheSchwarz})
are related to the module of second-order differential operators.

\begin{pro}
If $\dim M\geq2$,
the action of $\Diff(M)$ on the space of the space $\cD^2_\l$
of second-order differential operators, defined by
(\ref{Oaction},\ref{TheDiagram}) is as follows~:
\begin{eqnarray}
\left(\bar f_\l\bar a_2\right)^{ij}&=&
(f^*\bar a_2)^{ij}
\nonumber\\
\label{ExplAction}
\left(\bar f_\l\bar a_1\right)^i&=&
(f^*\bar a_1)^i\;+\;
\displaystyle
(2\l-1)\,\frac{n+1}{n+3}\,\ell_{kl}^i(f^{-1})(f^*\bar a_2)^{kl}\\
\displaystyle
\bar f_\l\bar a_0&=&
f^*\bar a_0
\;-\;
\displaystyle
\frac{2\l(\l-1)}{n+2}\,S_{kl}(f^{-1})(f^*\bar a_2)^{kl}\nonumber
\end{eqnarray}
where $f^*$ is the natural action of $f$ on the symmetric
contravariant tensor fields.
\label{MainAct}
\end{pro}

\medskip
\noindent
\textbf{Remark}.
In the one-dimensional case, the formula (\ref{ExplAction})
holds true, recall that $\ell(f)\equiv0$ and 
$S_{kl}(f^{-1})(f^*\bar a_2)^{kl}=S(f^{-1})f^*\bar a_2$
with the operator of multiplication by the
classical Schwarzian derivative in the right hand side
(cf. \cite{bo}).
This shows that the cocycle (\ref{TheSchwarz}) is, indeed,
its natural generalization.

Note also, that the formula (\ref{fu}) is a particular case
of (\ref{ExplAction}).

\subsection{Module of differential operator as a deformation}

The space of differential operators $\cD^2_\l$
as a module over the Lie algebra of vector fields $\Vect(M)$
was first studied in \cite{do}, it was shown
that this module can be naturally considered as a deformation
of the module of tensor fields on $M$.
Proposition \ref{MainAct} extends this result to the level
of the diffeomorphism group $\Diff(M)$.
The formula (\ref{ExplAction}) shows that the $\Diff(M)$-module
of second order differential operators on $M$
$\cD^2_\l$ is a \textit{nontrivial deformation} of the module of
tensor fields $\cT^2$ generated by the cocycles
(\ref{ell}) and (\ref{TheSchwarz}).

In the one-dimensional case, the $\Diff(S^1)$-modules
of differential operators and the related higher order analogues of
the Schwarzian derivative was studied in \cite{bo}.

\bigskip

{\it Acknowledgments}. It is a pleasure to acknowledge 
numerous fruitful discussions with Christian~Duval and his constant
interest in this work, we are also grateful to
Pierre~Lecomte for fruitful discussions.

\vskip 1cm




\end{document}